\documentclass[10pt,twoside]{article}
\usepackage{amsmath,amsthm}
\usepackage{amssymb,latexsym}
\usepackage{Latex-document}

\title{\bf Algebraic Structures on Valuations,\vskip -2mm
Their Properties and Applications\vskip 6mm}

\author{Semyon Alesker\vspace*{-0.5cm}\thanks{Department of Mathematics, Tel Aviv University, Ramat Aviv,
69978 Tel Aviv, Israel. E-mail: semyon@post.tau.ac.il}}

\date{\vspace{-8mm}}

\newcommand{\RR}{\mbox{\rm $~\vrule height6.5pt width0.5pt
depth0.3pt\!\!$R}}

\newcommand{\CC}{\mbox{\rm $~\vrule height6.5pt width0.5pt
depth0.3pt\!\!$C}}

\def\eps{\varepsilon}

\def\lam{\lambda}

\def\str{\longrightarrow}

\def\agrc{{}\!^{\CC} {\cal A}Gr}

\swapnumbers
\newtheorem{theorem}{Theorem}[subsection]

\newtheorem{proposition}[theorem]{Proposition}

\theoremstyle{definition}

\newtheorem{definition}[theorem]{Definition}

\renewcommand{\thesection}{\arabic{section}.}
\renewcommand{\thesubsection}{\thesection\arabic{subsection}}

\def\addsec{\addtocounter{subsection}{1} \setcounter{theorem}{0}}

\begin{document}

\maketitle

\thispagestyle{first} \setcounter{page}{757}

\begin{abstract}

\vskip 3mm

We describe various structures of algebraic nature on the space of
continuous valuations on convex sets, their properties (like
versions of Poincar\'e duality and hard Lefschetz theorem), and
their relations and applications to integral geometry.

\vskip 4.5mm



\end{abstract}

\vskip 12mm

\setcounter{section}{-1}

\section{Introduction}
\setcounter{subsection}{1}

\vskip-5mm \hspace{5mm}

 The theory of continuous valuations on
convex sets generalizes, in a sense, both the measure theory and
the theory of the Euler characteristic. Roughly speaking one
should think of a continuous valuation $\phi$ on a real linear
space $V$ as a finite additive measure on a class of compact nice
subsets of $V$ (say piecewise smooth submanifolds with corners)
which satisfy the following additional property (instead of the
usual sigma-additivity): the restriction of $\phi$ to the subclass
of convex compact domains with smooth boundary extends by
continuity to the class ${\cal K}(V)$ of {\itshape all} convex
compact subsets of $V$. Here the continuity is understood in the
sense of the Hausdorff metric on ${\cal K}(V)$. Remind that the
Hausdorff metric $d_H$ on ${\cal K}(V)$ depends on the choice of
the Euclidean metric on $V$ and it is defined as follows:
$d_H(A,B):=\inf\{ \eps >0|A\subset (B)_\eps \mbox{ and } B\subset
(A)_\eps\},$ where $(U)_\eps$ denotes the $\eps$-neighborhood of a
set $U$.  This condition of continuity turns out to be very strong
restriction and has a lot of consequences on purely algebraic
level. These properties will be discussed in this paper. The
simplest examples of such valuations are any smooth measure on $V$
and the Euler characteristic. Also it turns out that one of the
main tools used recently  in investigations of valuations is the
representation theory of real reductive groups and the
Beilinson-Bernstein theory of $D$-modules.

Now let us give the formal definition of valuation.
\begin{definition}
a) {\it A function $\phi :{\cal K}(V) \str \CC$ is called a
valuation if for any $K_1, \, K_2 \in {\cal K}(V)$ such that their
union is also convex one has}
$$\phi(K_1 \cup K_2)= \phi(K_1) +\phi(K_2) -\phi(K_1 \cap K_2).$$

b) {\it A valuation $\phi$ is called continuous if it is
continuous with respect the Hausdorff metric on ${\cal K}(V)$.}
\end{definition}
The linear space of all continuous valuations on $V$ will be
denoted by $CVal(V)$. It is a Fr\'echet space with the topology of
uniform convergence on compact subsets of ${\cal K}(V)$. In
Section 1 we discuss its dense subspace of polynomial smooth
valuations $(PVal(V))^{sm}$ (it has the topology of inductive
limit of Fr\'echet spaces). It turns out that this space has a
natural structure of associative commutative unital algebra (when
the unity is the Euler characteristic).
In Section 2 we discuss the space $Val(V)$ of translation
invariant continuous valuations. Its dense subspace
$(Val(V))^{sm}$ of so called smooth valuations is a subalgebra of
$(PVal(V))^{sm}$. It has a natural grading and satisfies a version
of Poincar\'e duality. This property follows from the
Irreducibility Theorem 2.1.3 which is by itself key result in the
investigation of valuations (see Subsection 2.1). Moreover
{\itshape even} smooth translation invariant continuous valuations
form a graded subalgebra of $(Val(V))^{sm}$ and satisfy a version
of the hard Lefschetz theorem (Subsection 2.2). This property
turns out to be closely related to the cosine transform problem in
the (Gelfand style) integral geometry solved recently in
\cite{alesker-bernstein}. These properties of valuations turn out
to be useful to obtain new explicit classification results on
valuations with additional invariance properties. The classical
Hadwiger theorem describes explicitly $SO(n)$- and
$O(n)$-invariant translation invariant continuous valuations on
the Euclidean space $\RR^n$. The new result is the classification
of unitarily invariant translation invariant continuous valuations
on the Hermitian space $\CC^n$ (Subsection 2.3). The main
application of the classification results on valuations is
integral geometric formulas. Using our classification we obtain
new results in (Chern style) integral geometry of real
submanifolds of complex spaces (Section 3).

\section{General continuous valuations}
\setcounter{subsection}{1} \setcounter{theorem}{0}

\vskip-5mm \hspace{5mm}

 In order to study general continuous
valuations let us remind the definition of {\itshape polynomial}
valuation introduced by Khovanskii and Pukhlikov
\cite{khovanskii-pukhlikov1}, \cite{khovanskii-pukhlikov2}.
\begin{definition}
{\it A valuation $\phi$ is called polynomial of degree $d$ if for
every $K\in {\cal K}(V)$ the function $x\mapsto \phi (K+x)$ is a
polynomial on $V$ of degree at most $d$.}
\end{definition}
Note that valuations polynomial of degree 0 are called {\itshape
translation invariant} valuations. Polynomial valuations have many
nice combinatorial-algebraic properties
(\cite{khovanskii-pukhlikov1}, \cite{khovanskii-pukhlikov2}). Also
in \cite{alesker-ann} the author have classified explicitly
rotation invariant polynomial continuous valuations on a Euclidean
space.

Let us denote the space of polynomial continuous valuations on $V$
by $PVal(V)$. One has
\begin{proposition}[\cite{alesker-mult}]
The space $PVal(V)$ of polynomial continuous valuations is dense
in the space of all continuous valuations $CVal(V)$.
\end{proposition}
The proof of this proposition is rather simple; it is a tricky use
of a form of the Peter-Weyl theorem (for the orthogonal group
$O(n)$), and in particular the convexity is not used in any
essential way.

Let us remind the basic definition of a smooth vector for a
representation of a Lie group. Let $\rho$ be a continuous
representation of a Lie group $G$ in a Fr\'echet space $F$. A
vector $\xi \in F$ is called $G$-smooth if the map $g\mapsto
\rho(g)\xi$ is infinitely differentiable map from $G$ to $F$. It
is well known the the subset $F^{sm}$ of smooth vectors is a
$G$-invariant linear subspace dense in $F$. Moreover it has a
natural topology of a Fr\'echet space (which is stronger than that
induced from $F$),  and the representation of $G$ is $F^{sm}$ is
continuous.

 We will especially be interested in polynomial
valuations which are $GL(V)$-smooth. This space will be denoted by
$(PVal(V))^{sm}$.

{\bf Example.} Let $\mu$ be a measure on $V$ with a polynomial
density with respect to the Lebesgue measure. Let $A\in {\cal
K}(V)$ be a strictly convex compact subset with smooth boundary.
Then
$$\phi(K):=\mu (K+A)$$ is a continuous polynomial smooth
valuation (here $K+A:=\{k+a|k\in K,\, a\in A\}$).

Let us denote by ${\cal G}(V)$ the linear space  of valuations on
$V$ which are finite linear combinations of valuations from the
previous example. It can be shown (using Irreducibility Theorem
2.1.3) that ${\cal G}(V)$ is dense in $(PVal(V))^{sm}$.
 Let $W$ be another linear real vector space. Let
us define the exterior product $\phi \boxtimes  \psi \in {\cal G}
(V\times W)$ of two valuations $\phi \in{\cal G}(V), \, \psi\in
{\cal G}(W)$. Let $\phi (K)=\sum_i \mu_i(K+A_i),\,
\psi(L)=\sum_j\nu_j(L+A_j)$. Define
$$(\phi\boxtimes  \psi) (M):=\sum_{i,j}(\mu_i \boxtimes \nu
_j)(M+(A_i\times \{0\})+(\{0\}\times B_j)),$$ where $\mu_i
\boxtimes \nu _j$ denotes the usual product measure.
\begin{proposition}[\cite{alesker-mult}]
For $\phi \in{\cal G}(V), \, \psi\in {\cal G}(W)$ their exterior
product $\phi\boxtimes  \psi \in {\cal G}(V\times W)$ is well
defined; it is bilinear  with respect to each argument. Moreover
$$(\phi \boxtimes \psi)\boxtimes \eta= \phi \boxtimes (\psi\boxtimes
\eta).$$
\end{proposition}
Now let us define a product on ${\cal G}(V)$. Let $\Delta:
V\hookrightarrow V\times V$ denote the diagonal imbedding. For
$\phi, \, \psi \in {\cal G}(V)$ let
$$ \phi \cdot \psi :=\Delta^* (\phi \boxtimes  \psi),$$
where $\Delta ^*$ denotes the restriction of a valuation on
$V\times V$ to the diagonal.
\begin{proposition}[\cite{alesker-mult}]
The above defined multiplication uniquely extends by continuity to
$(PVal(V))^{sm}$. Then $(PVal(V))^{sm}$ becomes an associative
commutative unital algebra where the unit is the Euler
characteristic $\chi$.
\end{proposition}

\section{Translation invariant continuous valuations}

\vskip-5mm \hspace{5mm}

For a linear finite dimensional real vector space $V$ let us denote by $Val(V)$ the space of {\itshape translation
invariant} continuous valuations on $V$. This is a Fr\'echet space with respect to the topology of uniform
convergence on compact subsets of ${\cal K}(V)$. In this section we will discuss properties of this space.

\subsection*{2.1. Irreducibility theorem and Poincar\'e duality} \addsec
\vskip-5mm \hspace{5mm}

It was shown by P. McMullen \cite{mcmullen-euler} that the space
$Val(V)$ of translation invariant continuous valuations on $V$ has
a natural grading given by the degree of homogeneity of
valuations. Let us formulate this more precisely.
\begin{definition}
{\it A valuation $\phi$ is called homogeneous of degree $k$ if for
every convex compact set $K$ and for every scalar $\lam >0$}
$$\phi(\lam K)=\lam ^k \phi (K).$$
\end{definition}
Let us denote by $Val _k(V)$ the space of translation invariant
continuous valuations homogeneous of degree $k$.
\begin{theorem}[McMullen \cite{mcmullen-euler}]
$$Val(V)=\bigoplus_{k=0}^{n} Val_k(V),$$
where $n=\dim V$.
\end{theorem}

Note in particular that the degree of homogeneity is an integer
between 0 and $n=\dim V$. It is known that $Val_0(V)$ is
one-dimensional and is spanned by the Euler characteristic $\chi$,
and $Val_n(V)$ is also one-dimensional and is spanned by a
Lebesgue measure \cite{hadwiger-book}. The space $Val_n(V)$ is
also denoted by $Dens(V)$ (the space of densities on $V$). One has
further decomposition with respect to parity:
$$Val_k(V)=Val_k^{ev}(V)\oplus Val_k^{odd}(V),$$
where $Val_k^{ev}(V)$ is the subspace of even valuations ($\phi$
is called even if $\phi(-K)=\phi(K)$ for every $K\in {\cal
K}(V)$), and $Val_k^{odd}(V)$ is the subspace of odd valuations
($\phi$ is called odd if $\phi (-K)=-\phi(K)$ for every $K\in
{\cal K}(V)$). The Irreducibility Theorem is as follows.
\begin{theorem}[\cite{alesker-gafa},\cite{alesker-adv}]
The natural representation of the group $GL(V)$ on each space
$Val_k^{ev}(V)$ and $Val_k^{odd}(V)$ is irreducible.
\end{theorem}
This theorem is the main tool in further investigations of
valuations and classification of them (see Subsection 2.3). This
immediately implies so called McMullen's conjecture
\cite{mcmullen-conj}.
 Its proof is heavily based on the use of the representation theory of
real reductive groups and the Beilinson-Bernstein theory of
D-modules. Another key tool in the proof of this result is the
Klain-Schneider characterization of simple translation invariant
continuous valuations \cite{klain}, \cite{schneider}.

By the results of Section 1 $(Val(V))^{sm}$ is a subalgebra of
$(PVal(V))^{sm}$. It is easy to see that the algebra structure is
compatible with the grading, namely
$$(Val_i(V))^{sm}\otimes (Val_j(V))^{sm}\str
(Val_{i+j}(V))^{sm}.$$ In particular we have
$$(Val_i(V))^{sm}\otimes (Val_{n-i}(V))^{sm}\str
Dens(V).$$ A version of the Poincar\'e duality theorem says that
this is a perfect pairing. More precisely
\begin{theorem}[\cite{alesker-mult}]
The induced map
$$(Val_i(V))^{sm}\str (Val_{n-i}(V)^*)^{sm}\otimes Dens(V)$$
is an isomorphism.
\end{theorem}
\subsection*{2.2. Even translation invariant continuous valuations} \addsec
\vskip-5mm \hspace{5mm}

Let us denote by $Val^{ev}(V)$ the subspace of {\itshape even}
translation invariant continuous valuations. Then clearly
$(Val^{ev}(V))^{sm}$ is a subalgebra of $(Val(V))^{sm}$. It turns
out that it satisfies a version of the hard Lefschetz theorem
which we are going to describe.

Let us fix on $V$ a scalar product. Let $D$ denote the unit ball
with respect to this product. Let us define an operator $\Lambda:
Val(V)\str Val(V).$ For a valuation $\phi \in Val(V)$ set
$$(\Lambda \phi)(K):=\frac{d}{d\eps}|_{\eps =0} \phi (K+\eps D).$$
(Note that by a result of P. McMullen \cite{mcmullen-euler} $\phi
(K+\eps D)$ is a polynomial in $\eps
>0$ of degree at most $n$.) It is easy to see that $\Lambda$
preserves the parity of valuations and decreases the degree of
homogeneity by 1. In particular
$$\Lambda: Val_k^{ev}(V)\str Val_{k-1}^{ev}(V).$$
The following result is a version of the hard Lefschetz theorem.
\begin{theorem}[\cite{alesker-univa}]
Let $k>n/2$. Then
$$\Lambda
^{2k-n}:(Val^{ev}_k(V))^{sm}\str(Val^{ev}_{n-k}(V))^{sm}$$ is an
isomorphism. In particular for $1\leq i\leq 2k-n$ the map
$$\Lambda
^{i}:(Val^{ev}_k(V))^{sm}\str(Val^{ev}_{k-i}(V))^{sm}$$ is
injective.
\end{theorem}
Note that the proof of this result is based on the solution of the
cosine transform problem due to J. Bernstein and the author
\cite{alesker-bernstein}, which is the problem from (Gelfand
style) integral geometry motivated by stochastic geometry and
going back to G. Matheron \cite{matheron}.

\subsection*{2.3. Valuations invariant under a group} \addsec
\vskip-5mm \hspace{5mm}

Let $G$ be a subgroup of $GL(V)$. Let us denote by $Val^G(V)$ the
space of $G$-invariant translations invariant continuous
valuations. From the results of \cite{alesker-adv} and
\cite{alesker-univa} follows the following result.
\begin{theorem}
Let $G$ be a compact subgroup of $GL(V)$ acting transitively on
the unit sphere. Then $Val^G(V)$ is a finite dimensional graded
subalgebra of $(Val(V))^{sm}$. It satisfies the Poincar\'e
duality, and if $-Id\in G$ it satisfies the hard Lefschetz
theorem.
\end{theorem}
It turns out that $Val^G(V)$ can be described explicitly (as a
vector space) for $G=SO(n),\, O(n),\mbox{ and }U(n)$. In the first
two cases it is the classical theorem of Hadwiger
\cite{hadwiger-book}, the last case is new (see
\cite{alesker-univa}). In order to state these results we have to
introduce first sufficiently many examples.

Let $\Omega$ be a compact domain in a Euclidean space $V$ with a
smooth boundary $\partial \Omega$. Let $n=\dim V$. For any point
$s\in \partial \Omega$ let $k_1(s), \dots, k_{n-1}(s)$ denote the
principal curvatures at $s$. For $0\leq i \leq n-1$ define
$$V_i(\Omega):=\frac{1}{n} {n-1\choose n-1-i}^{-1}
\int _{\partial \Omega}\{k_{j_1}, \dots,k_{j_{n-1-i}}\} d\sigma,
$$ where $\{k_{j_1}, \dots,k_{j_{n-1-i}}\}$ denotes the $(n-1-i)$-th
elementary symmetric polynomial in the principal curvatures,
$d\sigma$ is the measure induced on $\partial \Omega$ by the
Euclidean structure. It is well known that $V_i$ (uniquely)
extends by continuity in the Hausdorff metric to ${\cal K}(V)$.
Define also $V_n(\Omega):=vol(\Omega)$. Note that $V_0$ is
proportional to the Euler characteristic $\chi$. It is well known
that $V_0, V_1, \dots ,V_n$  belong to $Val^{O(n)}(V)$. It is easy
to see that $V_k$ is homogeneous of degree $k$. The famous result
of Hadwiger says
\begin{theorem}[Hadwiger, \cite{hadwiger-book}]
Let $V$ be $n$-dimensional Euclidean space. The valuations $V_0,
V_1, \dots ,V_n$ form a basis of $Val^{SO(n)}(V)(=Val^{O(n)}(V))$.
\end{theorem}

Now let us describe unitarily invariant valuations on a Hermitian
space. Let $W$ be a Hermitian space, i.e. a complex vector space
equipped with a Hermitian scalar product. Let $m:=\dim _{\CC} W$
(thus $\dim _{\RR}W=2m$). For every non-negative integers $p$ and
$k$ such that $2p\leq k\leq 2m$  let us introduce the following
valuations:
$$U_{k,p}(K)=\int_{E\in \agrc _{m-p}} V_{k-2p}(K\cap E) \cdot
dE.$$ Then $U_{k,p}\in Val_k^{U(m)}(W)$.
\begin{theorem}[\cite{alesker-univa}]
Let $W$ be a Hermitian vector space of complex dimension $m$. The
valuations $U_{k,p}$ with $0\leq p \leq \frac{min\{k, 2m-k\}}{2}$
form a basis of the space $Val_k^{U(m)}(W)$.
\end{theorem}
It turns out that the proof of this theorem is highly indirect,
and it uses everything known about even translation invariant
continuous valuations including the solution of McMullen's
conjecture, cosine transform, hard Lefschetz theorem for
valuations, and also results of Howe and Lee \cite{howe-lee} on
the structure of certain $GL_n(\RR)$-modules. Namely in order to
describe explicitly the (finite dimensional) space of unitarily
invariant valuations it is necessary to study the (infinite
dimensional) $GL_{\RR}(W)$-module $Val^{ev}(W)$.

Note that as algebra $Val^{SO(n)} (V)$ is isomorphic to
$\CC[x]/(x^{n+1}$. The algebra structure of $Val^{U(m)}(W)$ is not
yet computed.

\section{Applications to integral geometry}
\setcounter{subsection}{1} \setcounter{theorem}{0}

\vskip-5mm \hspace{5mm}

 In this section we state new results from (Chern style) integral
geometry of Hermitian spaces. They are obtained by the author in
\cite{alesker-univa} using the classification of unitarily
invariant valuations described in Subsection 2.3 of this paper.
They can be considered as a generalization of the classical
kinematic formulas due to Chern, Crofton, Santal\'o,  and others
(see e.g. \cite{chern1},
\cite{chern2},\cite{griffiths},\cite{klain-rota}, \cite{santalo}).

Let us remind first the principal kinematic formula following
Chern \cite{chern1}. Let $ISO(n)$ denote the group of affine
isometries of the Euclidean space $\RR^n$. Let $\Omega_1,\,
\Omega_2$ be compact domains with smooth boundary in $\RR^n$.
Assume also that $\Omega_1 \cap U(\Omega_2)$ has finitely many
components for all $U\in ISO(n)$.
\begin{theorem}[\cite{chern1}]
$$\int_{U\in ISO(n)}\chi(\Omega_1\cap U(\Omega_2))dU=
\sum_{k=0}^{n} \kappa_k V_k(\Omega_1) V_{n-k}(\Omega_2),$$ where
$\kappa_k$ are constants depending on $k$ and $n$ only which can
be written down explicitly.
\end{theorem}
For the explicit form of the constants $\kappa_k$ we refer to
\cite{chern1} or \cite{santalo}, Ch.15 \S 4.

Let us return back to the Hermitian situation. Let $IU(m)$ denote
the group of affine isometries of the Hermitian space $\CC^m$
preserving the complex structure (then $IU(m)$ is isomorphic to
$\CC^m \rtimes U(m)$). Let $\Omega_1,\, \Omega_2$ be compact
domains with smooth boundary in $\CC^m$ such that $\Omega_1 \cap
U(\Omega_2)$ has finitely many components for all $U\in IU(m)$.
The new result is
\begin{theorem}[\cite{alesker-univa}]
$$\int_{U\in IU(m)} \chi (\Omega_1\cap U(\Omega_2)) dU=
\sum_{k_1+k_2=2m}\sum_{p_1, p_2} \kappa(k_1,k_2,p_1,p_2)
U_{k_1,p_1}(\Omega_1) U_{k_2,p_2}(\Omega_2),$$ where the inner sum
runs over $0\leq p_i\leq k_i/2, \, i=1,2$, and
$\kappa(k_1,k_2,p_1,p_2)$ are certain constants depending on
$m,k_1,k_2,p_1,p_2$ only.
\end{theorem}

{\bf Remark.} We could compute explicitly the constants
$\kappa(k_1,k_2,p_1,p_2)$ only in $\CC^2$.

For more integral geometric formulas of this and other type for
real domains in $\CC^m$ we refer to \cite{alesker-univa}.

\end{document}